\input amstex

\documentstyle{amsppt}

%\TagsOnRight

%\tt

\hsize14cm

\vsize19cm

\magnification=1200

\baselineskip=14pt

\NoBlackBoxes

\def\I{{\text{\rm I}}}

\def\II{{\text{\rm I\!I}}}

\def\R{\Bbb{R}}

\def\H{\Bbb{S}^4(1)}

\topmatter

\rm

\title
Non-closed minimal hypersurfaces of $\H$ with identically
zero Gau\ss-Kronecker curvature
\endtitle

\author
Tsasa Lusala
\endauthor

\thanks 2000 {\it Mathematics Subject Classification}. 53 B 25, 53 C 40.
\endthanks

\thanks This work was partially supported by the DFG-project Si-163/7-2 and a
DFG-NSFC exchange project.
\endthanks

\abstract
We give a partial local description of minimal hypersurfaces $M^3$ with
identically zero Gau\ss-Kronecker curvature function in the unit $4$-sphere $\H$, without assumption on the compactness of $M^3$.

\bigskip

\endabstract

\endtopmatter

\document

\head {\S 1. Introduction}
\endhead
Let $x\colon M^3\longrightarrow\H\subset\R^5$ be a
hypersurface immersion of a connected and orientable $3$-dimensional
manifold $M^3$ of class $C^{\infty}$ into
$\H\subset\R^5$. Let $\lambda_1$, $\lambda_2$ and
$\lambda_3$ be the three principal curvature functions. The normalized
elementary symmetric curvature functions of the immersion $x$ are given by:
$$\align
H&:=\frac{1}{3}\big(\lambda_1+\lambda_2+\lambda_3\big),\cr
H_2&:=\frac{1}{3}\big(\lambda_1\lambda_2+\lambda_1\lambda_3+\lambda_2\lambda_3\big),\cr
K&:=\lambda_1\lambda_2\lambda_3.
\endalign $$

S. Almeida and F. Brito \cite{1} suggested to classify closed hypersurface immersions for which two of the three functions $H$, $H_2$, $K$ are constant. The paper \cite{3} gives a survey of results on closed hypersurfaces in\/ $\H$ with two constant curvature functions.

\medskip

\noindent
Particularly, the paper \cite{2} investigated closed minimal
hypersurfaces with constant Gau\ss-Kronecker curvature function,
corresponding to $H\equiv0$ and $K\equiv const$. There it is proved that
closed minimal hypersurfaces with constant Gau\ss-Kronecker curvature
$K\not=0$ are isoparametric, therefore closed minimal hypersurfaces
with constant Gau\ss-Kronecker curvature $K\not=0$ are classified. Brito conjectured that all hypersurfaces in $\H$ with $K\equiv const\not=0$ and $H\equiv const$ (or $H_2\equiv const$) must be isoparametric (personal communication).
If $K\equiv0$ on $M^3$, the following is well known: a closed minimal
hypersurface immersion in $\H$ with nowhere zero second fundamental form
is a boundary of a tube which is built over a non-degenerate
minimal $2$-dimensional surface immersion in $\H$ with geodesic
radius $\frac{\pi}{2}$. This nice result proves the existence of
non-isoparametric closed minimal hypersurfaces with $K\equiv0$ in
$\H$. But so far no explicit non-isoparametric example has been given.
In this paper we investigate local descriptions of minimal hypersurfaces (not
necessarily closed) in $\H$ with identically zero Gau\ss-Kronecker curvature
$K$, but with nowhere zero second fundamental form. In particular we present the following two explicit non-isoparametric examples:

\medskip

\subheading{Example 1.1}
The mapping
$$\align
&x_1\colon \R^3\longrightarrow\H\subset\R^5\cr
&x_1(u, v, z)=\frac{1}{\sqrt{1+z^2}}\Big(\cos(\sqrt{2}u)C_1+\sin(\sqrt{2}u)C_2+\cos(\sqrt{2}v)C_3+\sin(\sqrt{2}v)C_4+zC_5\Big),
\endalign
$$
where $C_1,C_2,C_3,C_4,C_5\in\R^5$ are constant orthogonal vectors in
$\R^5$ such that $$\frac{1}{2}=<C_1, C_1>=<C_2, C_2>=<C_3, C_3>=<C_4,
C_4> \quad \text{and} \quad <C_5,C_5>=1,$$
defines a minimal hypersurface immersion with zero Gau\ss-Kroneker
curvature. The principal curvature functions take the values $\lambda_1(u,v,z)=\sqrt{z^2+1}$,
$\lambda_2(z)=-\sqrt{z^2+1}$ and $\lambda_3(z)=0$; they depend only on $z$.

\bigskip

\subheading{Example 1.2}
Let $I\subset\R$ be an open interval, $0<c_1, c_2\in\R$ such that $c_2e^{2v}-1-c_1^2e^{4v}>0$ for all $v\in I$, and $g, h\colon I\longrightarrow\R$ two differentiable functions on $I$ which are linearly independent solutions of the second order differential equation
$$(c_2e^{2v}-1-c_1^2e^{4v})A_5''(v)+(1-c_1^2e^{4v})A_5'(v)+2A_5(v)=0,$$
and such that $g^2(v)+h^2(v)=1-\frac{e^{-2v}}{c_2}$, for all $v\in I$. The existence of such functions will be proved below, see Lemma 3.4 and Remark 3.5. Then the mapping
$$\align
&x_2\colon \R\times I\times\R\longrightarrow\H\subset\R^5,\cr
&x_2(u, v, z)=\frac{e^{-v}}{\sqrt{c_2(z^2+1)}}\Big(\cos(u)C_1+\sin(u)C_2\Big)+\frac{1}{\sqrt{z^2+1}}\Big(zC_3+g(v)C_4+h(v)C_5\Big),
\endalign
$$
where $C_1, C_2, C_3, C_4, C_5\in\R^5$ are constant orthonormal vectors,
defines a minimal hypersurface immersion in $\H$ with identically zero
Gau\ss-Kronecker curvature. The principal curvatures take the values
$$\lambda_1(v,z)=c_1e^{2v}\sqrt{z^2+1}\quad , 
\lambda_2(v,z)=-c_1e^{2v}\sqrt{z^2+1} \quad \text{and}  \quad
\lambda_3(v,z)=0,$$
thus they depend only on $v$ and $z$.

\medskip

\noindent
The principal curvatures of both examples do depend on at most two parameters. We prove the following local classification of such hypersurfaces.

\medskip

\noindent
{\bf Main result:}
Let $x\colon M^3\longrightarrow\H\subset\R^4$ be a minimal hypersuface
immersion (with nowhere zero second fundamental form) of a connected and
orientable $C^{\infty}$-manifold $M^3$ in $\H$ with identically zero Gau\ss-Kronecker
curvature. If one of the two nowhere zero principal curvature functions is constant along its associated principal curvature line, then there exist local
coordinates so that the immersion $x$ locally can be described by one
of the two non-isoparametric hypersurfaces $x_1$ and $x_2$ (see Example 1.1 and Example 1.2) above, or locally by Cartan's minimal isoparametric hypersurface with principal curvatures $\sqrt{3}$, $-\sqrt{3}$ and $0$.

\head{\S 2. Notations and integrability conditions}
\endhead

Let $x:M^3\longrightarrow\H\subset\R^5$ be an immersion of a connected, orientable $3$-dimensional $C^{\infty}$-manifold $M^3$ into the unit $4$-sphere $\H$. Denote by $y$ a unit normal vector field on $\H$ along the immersion $x$, by $< , >$ the canonical inner product of the Euclidean structure, and by $\bar{\nabla}$ the flat connection of $\R^5$. Referring to \cite{4} for details on geometry of submanifolds, recall that as immersion of codimension $2$ in $\R^5$ the structure equations (Gau{\ss} and Weingarten equations) for $x$ state:
$$\left\{\eqalign{
\bar{\nabla}_udx(v)&:=dx(\nabla_uv)+\II(u,v)y-\I(u, v)x,\cr
dy(u)&:=-dx(Su),\cr}\right. \quad \text{for all} \quad u,v\in TM^3 \tag2.1
$$
where $\I$ denotes the first fundamental form (induced metric) with
Levi-Civita connection $\nabla$, $\II$ defines the second fundamental
form and $S$ denotes the shape operator.

\medskip

The structure equations imply the following integrability conditions (Gau{\ss} formula and Codazzi equation) for any $u,v,w\in TM^3$:
$$\align
R(u, v)w&=\I(w, v)u-\I(w, u)v+\II(w, v)Su-\II(w, u)Sv, \tag2.2 \cr
\big(\nabla_uS\big)v&=\big(\nabla_vS\big)u, \tag2.3
\endalign
$$
where $R$ denotes the Riemannian curvature tensor for the induced metric $\I$.

\medskip

Let $(e_1, e_2, e_3)$ be a  $\I$-orthonormal local differentiable
frame of principal curvature vector fields on $M^3$:
$$Se_1=\lambda_1 e_1, \quad Se_2=\lambda_2 e_2 \quad \text{and} \quad
Se_3=\lambda_3 e_3.$$
There are $9$ functions $\alpha_1,\cdot\cdot\cdot,\alpha_9$ such that
$$\left\{\eqalign{
\nabla_{e_1}e_1&=\alpha_1 e_2+\alpha_2 e_3, \quad \nabla_{e_1}e_2=-\alpha_1 e_1+\alpha_3 e_3, \quad
\nabla_{e_1}e_3=-\alpha_2 e_1-\alpha_3 e_2;\cr
\nabla_{e_2}e_1&=-\alpha_4 e_2+\alpha_6 e_3, \quad
\nabla_{e_2}e_2=\alpha_4 e_1+\alpha_5 e_3, \quad
\nabla_{e_2}e_3=-\alpha_6 e_1-\alpha_5 e_2;\cr
\nabla_{e_3}e_1&=\alpha_9 e_2-\alpha_7e_3, \quad
\nabla_{e_3}e_2=-\alpha_9 e_1-\alpha_8e_3, \quad
\nabla_{e_3}e_3=\alpha_7 e_1+\alpha_8e_2.\cr}\right. \tag2.4
$$

\medskip
\subheading{Remark 2.1}
Consider the situation that the three principal curvature functions
$\lambda_1,\lambda_2,\lambda_3$ are everywhere distinct; the fact that
the frame $(e_1,e_2,e_3)$ is orthonormal implies that the functions
$\alpha_i, i=1,\cdot\cdot\cdot,9$ in (2.4) are defined on $M^3$ uniquely up to sign.

Applying the Codazzi equation (2.3) to the vector fields $e_1,
e_2, e_3$ and using (2.4), one gets the following equations:
$$\left\{\eqalign{
e_1(\lambda_2)&=\alpha_4(\lambda_2-\lambda_1), \cr
e_1(\lambda_3)&=\alpha_7(\lambda_3-\lambda_1), \cr
e_2(\lambda_1)&=\alpha_1(\lambda_1-\lambda_2), \cr
e_2(\lambda_3)&=\alpha_8(\lambda_3-\lambda_2), \cr
e_3(\lambda_1)&=\alpha_2(\lambda_1-\lambda_3), \cr
e_3(\lambda_2)&=\alpha_5(\lambda_2-\lambda_3),\cr
\alpha_9(\lambda_1-\lambda_2)&=\alpha_3(\lambda_2-\lambda_3)=\alpha_6(\lambda_1-\lambda_3).\cr}\right. \tag2.5
$$

From now we assume that the immersion $x$ (with nowhere zero second
fundamental form) is minimal and has identically zero Gau\ss-Kronecker curvature ($K\equiv0$). There exists a positive non-zero function $\lambda$ such that the principal curvature functions associated to the immersion $x$ are $\lambda_1=\lambda$, $\lambda_2=-\lambda$ and $\lambda_3=0$. From the equations (2.5), one gets
$$\align
e_1(\lambda)&=2\alpha_4\lambda, \quad \quad e_2(\lambda)=2\alpha_1\lambda, \quad \quad e_3(\lambda)=\alpha_2\lambda;\tag2.6\cr
\alpha_5&=\alpha_2, \quad \quad 2\alpha_9=-\alpha_3=\alpha_6, \quad \quad \alpha_7=0=\alpha_8. \tag2.7
\endalign
$$
Applying the Gau{\ss} formula (2.2) to the vector fields $e_1,e_2,e_3$
and using the equations (2.4) and (2.7), one gets
$$\left\{\eqalign{
e_1(\alpha_4)+e_2(\alpha_1)&=1-\lambda^2+\alpha_1^2+\alpha_2^2+2\alpha_3^2+\alpha_4^2  \cr
e_3(\alpha_1)+\frac{1}{2}e_1(\alpha_3)&=\alpha_1\alpha_2-\frac{1}{2}\alpha_3\alpha_4 \cr
e_3(\alpha_4)-\frac{1}{2}e_2(\alpha_3)&=\alpha_2\alpha_4+\frac{1}{2}\alpha_1\alpha_3, \cr
e_3(\alpha_2)&=1+\alpha_2^2-\alpha_3^2, \cr
e_3(\alpha_3)&=2\alpha_2\alpha_3, \cr
e_1(\alpha_2)&=e_2(\alpha_3), \cr
e_1(\alpha_3)&=-e_2(\alpha_2). \cr} \right.\tag2.8
$$

Note that the Lie brackets with respect to the vector fields $e_1$, $e_2$ and $e_3$ are given by:
$$
[e_1,e_2]=-\alpha_1e_1+\alpha_4e_2+2\alpha_3e_3, \quad [e_1,e_3]=-\alpha_2e_1-\frac{1}{2}\alpha_3e_2, \quad [e_2,e_3]=\frac{1}{2}\alpha_3e_1-\alpha_2e_2.
$$

The fundamental equations (2.1) applied to the vector fields $e_1, e_2, e_3$ give rise to the following (partial) differential equations:

$$
\left\{\eqalign{\bar{\nabla}_{e_1}dx(e_1)&=\alpha_1 dx(e_2)+\alpha_2dx(e_3)+\lambda y-x\cr
\bar{\nabla}_{e_1}dx(e_2)&=-\alpha_1 dx(e_1)+\alpha_3dx(e_3)\cr
\bar{\nabla}_{e_1}dx(e_3)&=-\alpha_2 dx(e_1)-\alpha_3dx(e_2)\cr
\bar{\nabla}_{e_2}dx(e_1)&=-\alpha_4 dx(e_2)-\alpha_3dx(e_3)\cr
\bar{\nabla}_{e_2}dx(e_2)&=\alpha_4 dx(e_1)+\alpha_2dx(e_3)-\lambda y-x\cr
\bar{\nabla}_{e_2}dx(e_3)&=\alpha_3 dx(e_2)-\alpha_2dx(e_2)\cr
\bar{\nabla}_{e_3}dx(e_1)&=-\frac{1}{2}\alpha_3 dx(e_2)\cr
\bar{\nabla}_{e_3}dx(e_2)&=\frac{1}{2}\alpha_3 dx(e_1)\cr
\bar{\nabla}_{e_3}dx(e_3)&=-x\cr
dy(e_1)&=-\lambda dx(e_1),\cr
dy(e_2)&=\lambda dx(e_2), \cr
dy(e_3)&=0.\cr}\right.\tag2.9
$$

\head {\S 3. Proof of the main result}
\endhead

To describe locally hypersurface immersions in \/ $\H$ with
$K\equiv0$, one has to find local coordinates to solve the structure
equations (2.9) using the integrability conditions
(2.6) and (2.8). It seems to be very difficult to
solve this problem in full generality.

\medskip

\noindent
In this section, we consider natural additional assumptions on the
functions $\alpha_1$, $\alpha_2$, $\alpha_3$ and $\alpha_4$ to solve
the structure equations (fundamental equations) for minimal
hypersurface immersions in $\H$ with $K\equiv0$; namely we assume
that the function $\lambda$ is constant along the $e_1$-direction. This
additional assumption is suggested by the Examples 1.1 and 1.2.
\subheading{Proposition 3.1}
Let $x\colon M^3\longrightarrow\H$ be a minimal hypersurface immersion
with identically zero Gau\ss-Kronecker curvature function and nowhere zero
second fundamental form. Let $\alpha_1,\cdot\cdot\cdot,\alpha_9$ be the the functions as defined in (2.4).
If the function $\alpha_4$ vanishes identically on an open subset $U$
of $M^3$, i.e., the function $\lambda$ is constant along the curvature
line of the vector field $e_1$, then also the function $\alpha_3$ vanishes identically on $U$, or the immersion $x$ is a minimal Cartan isoparametric hypersurface on $U$.

\subheading{\bf Proof}
Using the equations (2.6) and
$$e_1e_2(\lambda)-e_2e_1(\lambda)-[e_1, e_2](\lambda)=0=e_1e_3(\lambda)-e_3e_1(\lambda)-[e_1, e_3](\lambda),$$
 one gets:
$$e_1(\alpha_1)=\alpha_2\alpha_3\quad \text{and}\quad e_1(\alpha_2)=-\alpha_1\alpha_3.$$
Similarly, using (2.8) and $$e_2e_3(\alpha_3)-e_3e_2(\alpha_3)-[e_2, e_3](\alpha_3)=0,$$ one has:
$$\alpha_3e_2(\alpha_2)=0.$$
Assume now that $\alpha_3\not=0$ everywhere. This implies that
$e_2(\alpha_2)$ vanishes identically. Inserting again the equations (2.6) into
$$e_1e_2(\alpha_3)-e_2e_1(\alpha_3)-[e_1, e_2](\alpha_3)=0,$$
one gets:
$$\alpha_3^2\alpha_2=0.$$
Therefore $\alpha_2=0=\alpha_1$ and $\alpha_3=\pm1$. Consequently, $\lambda^2=3$. Thus the immersion is isoparametric with principal curvatures $\lambda_1=\sqrt{3}$,  $\lambda_2=-\sqrt{3}$ and $\lambda_3=0$, i.e. the immersion is a Cartan's minimal isoparametric hypersurface in $\H$.

\subheading{Corollary 3.2}
Let $x\colon M^3\longrightarrow\H\subset\R^5$ be a closed minimal hypersurface immersion of a connected and orientable
manifold $M^3$ into $\H\subset\R^5$ with nowhere zero second fundamental
form and $K\equiv0$. Assume that one of the functions $\alpha_1$ and $\alpha_4$ vanish identically on $M^3$. Then $x(M)$ is a Cartan's minimal isoparametric hypersurface of\/ $\H$, i.e., the boundary of the tube Tube($V^2, \frac{\pi}{2}$) with radius $\frac{\pi}{2}$ around the Veronese surface $V^2\subset\H$.

\subheading{Proof}
From the proposition above we have two possibilities:

\noindent

\itemitem{(i)} $\alpha_3$ vanishes identically on $M^3$;

\noindent

\itemitem{(ii)} or $\alpha_1=\alpha_2=\alpha_4\equiv 0$ and $\alpha_3^2=1$.

\medskip

\noindent
Assuming that the hypersurface $M^3$ is closed, the case (i) above cannot
happen because the function $e_3(\alpha_2)=1+\alpha_2^2-\alpha_3^2$ should be zero at the minimum and maximum points of the function $\alpha_2$; but with $\alpha_3\equiv0$,  $e_3(\alpha_2)=1+\alpha_2^2$ has no zeros.
\hfill$\square$

\subheading{Proposition 3.3}
Let $x\colon M^3\longrightarrow\H\subset\R^5$ be a minimal hypersurface immersion of a connected and orientable
manifold $M^3$ into $\H\subset\R^5$ with nowhere zero second fundamental
form and $K\equiv0$. Assume that
the immersion is non-isoparametric and the functions $\alpha_1$ and
$\alpha_4$ vanish identically on $M^3$. Then there are local
coordinates so that the immersion $x$ can be locally described by the
parametrization of the hypersurface given in Example 1.1.

\subheading{Proof}
From Proposition 3.1 we may assume that the function $\alpha_3$ vanishes
identically on $M^3$. Then the following equations hold:

$$\align
e_2(\alpha_2)&=0=e_1(\alpha_2)=e_1(\lambda)=e_2(\lambda),\cr
e_3(\alpha_2)&=\alpha_2^2+1,\cr 
e_3(\lambda)&=\lambda\alpha_2, \cr
\lambda^2&=\alpha_2^2+1.
\endalign
$$
The vector fields $\frac{1}{\lambda}e_1$, $\frac{1}{\lambda}e_2$ and $\frac{1}{\alpha_2^2+1}e_3$ satisfy:
$$0=\Big[\frac{1}{\lambda}e_1,\frac{1}{\lambda}e_2\Big]=\Big[\frac{1}{\lambda}e_1, \frac{1}{\alpha_2^2+1}e_3\Big]=\Big[\frac{1}{\lambda}e_2, \frac{1}{\alpha_2^2+1}e_3\Big].$$
Therefore there are local coordinates $(u,v,z)$ on $M^3$ such that
$$\frac{1}{\lambda}e_1=\frac{\partial}{\partial u}, \quad \frac{1}{\lambda}e_2=\frac{\partial}{\partial v}, \quad \frac{1}{\alpha_2^2+1}e_3=\frac{\partial}{\partial z}.$$
The foregoing equations give
$$\alpha_2=z \quad \text{and} \quad \lambda=\sqrt{z^2+1}.$$
With respect to the frame $(\frac{\partial}{\partial u},
\frac{\partial}{\partial v}, \frac{\partial}{\partial z})$ the structure
equations (Gau{\ss} and Weingarten equations) are given by the
following system of second order partial differential equations:

$$\align
\lambda^2x_{uu}&=z(z^2+1)x_z+\lambda y-x, \tag3.10\cr
x_{vu}&=0=x_{uv}, \tag3.11\cr
x_{zu}&=\frac{-z}{z^2+1}x_u=x_{uz}, \tag3.12\cr
\lambda^2x_{vv}&=z(z^2+1)x_z-\lambda y-x, \tag3.13\cr
x_{vz}&=\frac{-z}{z^2+1}x_v=x_{vz}, \tag3.14\cr
(z^2+1)^2x_{zz}&=-2z(z^2+1)x_z-x, \tag3.15\cr
y_u&=-\lambda x_u, \tag3.16\cr
y_v&=\lambda x_v, \tag3.17\cr
y_z&=0. \tag3.18
\endalign
$$

Differentiating (3.10) with respect to $u$ and using (3.12) and (3.16), one gets
$$x_{uuu}=-2x_u.$$
There are two vector valued functions $A_1\equiv A_1(v,z)$ and  $A_2\equiv A_2(v,z)$ in $\R^5$ depending only on $v$ and $z$ such that
$$
x_u=\sqrt{2}\Big(-\sin(\sqrt{2}u)A_1+\cos(\sqrt{2}u)A_2\Big).\tag3.19
$$
One has
$$\align
\frac{1}{\lambda^2}&=\I(x_u, x_u)\cr
&=2\Big(<A_1, A_1>\sin^2(\sqrt{2}u)-<A_1, A_2>\sin(2\sqrt{2}u)+<A_2, A_2>\cos^2(\sqrt{2}u)\Big)\cr
&=<A_1, A_1>+<A_2, A_2>+(<A_2, A_2>-<A_1, A_1>)\cos(2\sqrt{2}u)\cr
&\quad -2<A_1, A_2>\sin(2\sqrt{2}u).
\endalign
$$
The linear independence of the functions $1$, $\sin(2\sqrt{2}u)$ and
$\cos(2\sqrt{2}u)$ implies

$$
<A_1, A_1>=\frac{1}{2\lambda^2}=<A_2, A_2> \quad \text{and} \quad <A_1, A_2>=0.
\tag3.20$$
Furthermore there is a vector valued function $A_3\equiv A_3(v,z)$ depending
only on $v$ and $z$ such that $$x(u, v, z)=\cos(\sqrt{2}u)A_1(v, z)+\sin(\sqrt{2}u)A_2(v, z)+A_3(v, z).$$
One has
$$\align
1&=<x, x>\cr
&=\frac{1}{2\lambda^2}+2<A_1, A_3>\cos(\sqrt{2}u)+2<A_2, A_3>\sin(\sqrt{2}u)+<A_3, A_3>.
\endalign$$
From the linear independence of the functions $1$, $\sin(\sqrt{2}u)$ and $\cos(\sqrt{2}u)$, one gets
$$\align
1=\frac{1}{2\lambda^2}+<A_3, A_3>\quad \text{and} \quad <A_1, A_3>=0=<A_2, A_3>.\tag3.21
\endalign$$
Differentiating (3.19) with respect to $z$ (and with respect
to $v$, resp.) and using the equation (3.12) (the equation
(3.11), resp.) and the linear independence of the functions
$1$, $\sin(\sqrt{2}u)$ and $\cos(\sqrt{2}u)$, one gets the following
first order partial differential equations for the vector valued functions
$A_1(v, z)$ and $A_2(v, z)$:
$$
(z^2+1)\frac{\partial A_1}{\partial z}=-zA_1,  \quad \quad
(z^2+1)\frac{\partial A_2}{\partial z}=-zA_2 \quad \text{and} \quad
\frac{\partial A_1}{\partial v}=0=\frac{\partial A_2}{\partial v}.
$$
There are constant vectors $C_1$ and $C_2$ in $\R^5$ such that
$$
A_1=\frac{C_1}{\sqrt{z^2+1}} \quad \text{and} \quad A_2=\frac{C_2}{\sqrt{z^2+1}}.$$
Because of (3.20), one has
$$<C_1,C_1>=\frac{1}{2}=<C_2,C_2> \quad \text{and} \quad <C_1,C_2>=0.$$
The immersion $x$ takes the form
$$x(u, v, z)=\frac{1}{\sqrt{z^2+1}}\Big(\cos(\sqrt{2}u)C_1+\sin(\sqrt{2}u)C_2\Big)+A_3(v, z).$$
After inserting the above expression for $x$ into the equation
(3.15), one gets the following second order partial
differential equation for the vector valued function $A_3(v, z)$:
$$
(z^2+1)^2\frac{\partial^2 A_3}{\partial z^2}+2z(z^2+1)\frac{\partial A_3}{\partial z}+A_3=0.$$
There are two vector valued functions $A_4\equiv A_4(v)$ and $A_5\equiv
A_5(v)$ in $\R^5$ such that
$$A_3(v,z)=\frac{A_4(v)}{\sqrt{z^2+1}}+\frac{zA_5(v)}{\sqrt{z^2+1}}.$$
The equation (3.14) implies that the vector valued function $A_5(v)$
is constant: $$A_5(v)=C_5\equiv const.$$ Because of (3.21), one has
$$0=<C_5, A_4(v)>=<C_1, A_4(v)>=<C_2, A_4(v)>=<C_1, C_5>=<C_2, C_5>,$$
$$<A_4(v), A_4(v)>=\frac{1}{2} \quad \text{and} \quad <C_5, C_5>=1.$$

Now eliminating $y$ from the equations (3.10) and (3.13), one has that the vector valued function $A_4(v)$ is a solution of the following second order differential equation:
$$A_4''(v)=-2A_4(v).$$
Therefore there are constant vectors $C_3, C_4\in\R^3$ such that
$$A_4(v)=\cos(\sqrt{2}v)C_3+\sin(\sqrt{2}v)C_4.$$
The constant vectors $C_3, C_4$ are orthogonal to $C_1,C_2,C_5$ and satisfy
$$<C_3,C_3>=\frac{1}{2}=<C_4,C_4> \quad \text{and} \quad <C_3,C_4>=0.$$
Finally, the local description of the immersion $x$ is given by
$$x(u, v, z)=\frac{1}{\sqrt{z^2+1}}\Big(\cos(\sqrt{2}u)C_1+\sin(\sqrt{2}u)C_2+\cos(\sqrt{2}v)C_3+\sin(\sqrt{2}v)C_4+zC_5\Big).\tag3.22$$

Note that the vector field $y$ defined by
$$y(u, v,
z)=-\cos(\sqrt{2}u)C_1-\sin(\sqrt{2}u)C_2+\cos(\sqrt{2}v)C_3+\sin(\sqrt{2}v)C_4
$$
is unit and normal to $x$. It is then easy to check that the mapping
(3.22) defines a minimal hypersurface immersion with $K\equiv0$ in
$\H$. Clearly this hypersurface immersion is non-isoparametric. \hfill$\square$

\bigskip

Now we want to characterize Example 1.2. In order to succeed, we need to prove the following lemma:
\subheading{Lemma 3.4}
Let $I\subset\R$ be an open interval and $0<c_1, c_2\in\R$ be two constant
positive real numbers such that $c_2e^{2v}-1-c_1^2e^{4v}>0$ for every $v\in
I$. Consider the following second order differential equation for some function $A$ on $I$:
$$(c_2e^{2v}-1-c_1^2e^{4v})A''(v)+(1-c_1^2e^{4v})A'(v)+2A(v)=0.\tag3.23$$
Then there exists a function $\phi\colon I\longrightarrow\R$ such that the general solution $A(v)$ for the equation (3.23) takes the form
$$A(v)=a\sqrt{c_2-e^{-2v}}\cos\big(\phi(v)\big)+b\sqrt{c_2-e^{-2v}}\sin\big(\phi(v)\big),$$
where $a, b\in\R$ are constants.

\subheading{Proof}
Let $B\colon I\longrightarrow\R$ be the function defined by
$$A(v)=\colon\sqrt{c_2-e^{-2v}}B(v).$$ Inserting this expression of $A(v)$ into the equation (3.23), one gets the following second order differential equation in $B(v)$:
$$B''(v)-\frac{(3+c_1^2e^{4v}+c_2c_1^2e^{6v}-3c_2e^{2v})B'(v)}{(c_2e^{2v}-1)(c_2e^{2v}-1-c_1^2e^{4v})}+\frac{c_2c_1^2e^{6v}B(v)}{(c_2e^{2v}-1)^2(c_2e^{2v}-1-c_1^2e^{4v})}=0.\tag3.24$$
Now take $\phi$ to be the function on $I$ such that $$\phi'(v)=\sqrt{\frac{c_2c_1^2e^{6v}}{(c_2e^{2v}-1)^2(c_2e^{2v}-1-c_1^2e^{4v})}}.$$
It follows that
$$\frac{\phi''(v)}{\phi'(v)}=\frac{(3+c_1^2e^{4v}+c_2c_1^2e^{6v}-3c_2e^{2v})}{(c_2e^{2v}-1)(c_2e^{2v}-1-c_1^2e^{4v})}.$$
Thus the equation (3.24) becomes
$$B''(v)-\frac{\phi''(v)}{\phi'(v)}B'(v)+{\phi'}^2(v)B(v)=0.$$
Therefore, there are constants $a,b\in\R$ such that
$$B(v)=a\cos\big(\phi(v)\big)+b\sin\big(\phi(v)\big).$$ \hfill$\square$

\subheading{Remark 3.5}
The particular solutions $g(v)$ and $h(v)$ of the equation (3.23) given by
$$g(v)=\frac{1}{\sqrt{c_2}}\sqrt{c_2-e^{-2v}}\cos\big(\phi(v)\big) \quad \text{and} \quad g(v)=\frac{1}{\sqrt{c_2}}\sqrt{c_2-e^{-2v}}\sin\big(\phi(v)\big)$$
are linearly independent and satisfy $g^2(v)+h^2(v)=1-\frac{e^{-2v}}{c_2}$.

\subheading{Proposition 3.6}
Let $x\colon M^3\longrightarrow\H\subset\R^5$ be a minimal hypersurface 
immersion of a connected and orientable manifold $M^3$ into
$\H\subset\R^5$ with identically zero Gau\ss-Kronecker curvature, but with nowhere zero second fundamental form. Assume
that the function $\alpha_4$ vanishes identically and $\alpha_1$ is nowhere zero on $M^3$. Then there are local coordinates $(u, v, z)$ such that the immersion $x$ can be locally described by the parametrization of the hypersurface given in Example 1.2.

\subheading{Proof}
The vector fields $\frac{1}{\alpha_1}e_2$ and $\frac{1}{\alpha_2^2+1}e_3$ satisfy
$$[\frac{1}{\alpha_1}e_2, \frac{1}{\alpha_2^2+1}e_3]=0.$$
Define the function $f$ on $M^3$ by
$$
f:=\frac{1}{\sqrt{1+\lambda^2+\alpha_1^2+\alpha_2^2}}.\tag3.25
$$
The function $f$ satisfies the following equations:
$$
e_1(f)=0, \quad e_2(f)=-\alpha_1 f \quad \text{and} \quad e_3(f)=-\alpha_2 f.$$
Consequently,
$$
[f e_1, \frac{1}{\alpha_1}e_2]=0=[f e_1, \frac{1}{\alpha_2^2+1}e_3].$$
Therefore there are local coordinates $(u, v, z)$ on $M^3$ such that
$$\frac{\partial }{\partial z}=\frac{1}{\alpha_2^2+1}e_3, \quad
\frac{\partial }{\partial v}=\frac{1}{\alpha_1}e_2, \quad
\frac{\partial }{\partial u}=f e_1.$$
From above we get the following equations for $\alpha_2$, $\lambda$
and $\alpha_1$.

\medskip

\noindent
Equations for $\alpha_2$:
$$
\frac{\partial \alpha_2}{\partial u}=fe_1(\alpha_2)=0, \quad \frac{\partial \alpha_2}{\partial v}=\frac{1}{\alpha_1}e_2(\alpha_2)=0, \quad \frac{\partial \alpha_2}{\partial z}=\frac{1}{1+\alpha_2^2}e_3(\alpha_2)=1.
$$
So
$$
\alpha_2=z.\tag3.26$$
Equations for $\lambda$:
$$\frac{\partial \lambda}{\partial u}=0, \quad \frac{\partial \lambda}{\partial v}=2\lambda, \quad \text{and} \quad \frac{\partial \lambda}{\partial z}=\frac{z\lambda}{z^2+1}.$$
Therefore
$$\lambda\equiv\lambda(v,z)=c_1 e^{2v}\sqrt{z^2+1},\tag3.27$$
where $0<c_1\in\R$.

\medskip

\noindent
Equations for $\alpha_1$:
$$
\frac{\partial \alpha_1}{\partial z}=\frac{z\alpha_1}{z^2+1}, \quad \frac{\partial \alpha_1}{\partial v}=\frac{\alpha_1^2+z^2+1-\lambda^2}{\alpha_1}, \quad \text{and} \quad \frac{\partial \alpha_1}{\partial u}=0;$$
thus
$$\alpha_1^2\equiv\alpha_1^2(v,z)=(z^2+1)(c_2 e^{2v}-1-c_1^2e^{4v}),\tag3.28$$
where $v\in I\subset\R$ an open interval and $0<c_2\in\R$ is a constant such that $c_2 e^{2v}-1-c_1^2e^{4v}>0$ for all $v\in I$. Inserting the expressions (3.26), (3.27) and (3.28) (of $\alpha_2$, $\lambda$ and $\alpha_1$, respectively) into (3.25), we see that the function $f$ satisfies
$$(1+z^2)f^2(v)=\frac{e^{-2v}}{c_2}.$$

\bigskip

\noindent
With respect to the frame $(\frac{\partial}{\partial u}, \frac{\partial}{\partial v}, \frac{\partial}{\partial z})$ the structure equations are given by the following system of second order partial differential equations:
$$\align
x_{vz}&=-\frac{\alpha_2}{\alpha_2^2+1} x_v=x_{zv}, \tag3.29\cr
x_{zu}&=-\frac{\alpha_2}{\alpha_2^2+1} x_u=x_{uz}, \tag3.30\cr
x_{vu}&=-x_u=x_{uv}, \tag 3.31\cr
x_{uu}&=f^2\left(\alpha_1^2 x_v+\alpha_2(\alpha_2^2+1)x_z+\lambda y-x\right),
\tag3.32 \cr
x_{zz}&=-\frac{2\alpha_2}{\alpha_2^2+1}x_z-\frac{1}{(\alpha_2^2+1)^2}x, \tag3.33 \cr
x_{vv}&=\frac{1}{\alpha_1^2}\Big(-(\alpha_1^2+\alpha_2^2+1-\lambda^2)x_v+\alpha_2(\alpha_2^2+1)x_z-\lambda y-x\Big), \tag3.34 \cr
y_u&=-\lambda x_u,  \tag3.35\cr
y_v&=\lambda x_v, \tag3.36 \cr
y_z&=0. \tag3.37
\endalign
$$
Differentiating the equation (3.32) with respect to $u$, one has:
$$x_{uuu}=f^2(-\alpha_1^2-\alpha_2^2-\lambda^2-1)x_u=-x_u.\tag3.38 $$
 There are vector valued functions $A_1\equiv A_1(v, z)$ and $A_2\equiv A_2(v, z)$ in $\R^5$ depending only on $v$ and $z$ such that
$$x_u=-\sin(u)A_1+\cos(u)A_2.\tag3.39$$
One has:
$$\align
f^2&=\I(x_u, x_u)\cr
&=<A_1,A_1>\sin^2(u)-<A_1, A_2>\sin(2u)+<A_2,A_2>\cos^2(u)\cr
&=\frac{<A_2,A_2>+<A_1,A_1>}{2}+\frac{<A_2,A_2>-<A_1,A_1>}{2}\cos(2u)\cr
&\quad-<A_1, A_2>\sin(2u).
\endalign$$
Using the linear independence of the functions $1$, $\cos(2u)$ and
$\sin(2u)$, from the equation above we deduce
$$<A_1,A_1>=f^2=<A_2,A_2> \quad \text{and} \quad <A_1,A_2>=0. \tag3.40$$
Furthermore, there is a vector valued function $A_3\equiv A_3(v, z)$ depending only on $v$ and $z$ such that
$$x(u, v, z)=\cos(u)A_1(v, z)+\sin(u)A_2(v, z)+A_3(v, z). \tag3.41$$
One has
$$\align
1&=<x,x>\cr
&=f^2+<A_3,A_3>+2<A_1, A_3>\cos(u)+2<A_2, A_3>\sin(u).
\endalign$$
Using the linear independence of the functions $1$, $\cos(u)$ and $\sin(u)$, we get
$$
1=f^2+<A_3,A_3> \quad \text{and} \quad <A_1, A_3>=0=<A_2, A_3>. \tag3.42$$

\medskip

\noindent
We differentiate the equation (3.39) with respect to $z$ and $v$ and use the equations (3.30) and (3.31); then the linear independence of the functions $\cos(u)$ and $\sin(u)$ implies the
following first order partial differential equations for the vector valued
functions $A_1(v, z)$ and $A_2(v, z)$:
$$\align
\frac{\partial A_1}{\partial z}&=-\frac{z}{z^2+1} A_1,  \quad \frac{\partial
A_1}{\partial v}=-A_1,\tag3.43 \cr
\frac{\partial A_2}{\partial z}&=-\frac{z}{z^2+1} A_2,  \quad \frac{\partial
A_2}{\partial v}=-A_2.\tag3.44
\endalign$$
Therefore,
$$A_1(v, z)=f(v, z)C_1\quad \text{and} \quad A_2(v, z)=f(v, z)C_2, \tag3.45$$
where $C_1, C_2\in\R^5$ are constant vectors; they are orthonormal because of (3.40), and from (3.42) they are orthogonal to $A_3$.

\medskip

\noindent
Differentiating (3.41) with respect to $z$ we have:
$$\align
x_z&=\frac{\partial f}{\partial z}\Big(\cos(u)C_1+\sin(u)C_2\Big)+\frac{\partial A_3}{\partial z},\cr
x_{zz}&=\frac{\partial^2 f}{\partial z^2}\Big(\cos(u)C_1+\sin(u)C_2\Big)+\frac{\partial^2 A_3}{\partial z^2}.
\endalign$$
Using the equation (3.33), we get the following partial differential equation for $A_3$, depending only on $z$:
$$\frac{\partial^2 A_3}{\partial z^2}=-\frac{2z}{z^2+1}\frac{\partial A_3}{\partial z}-\frac{1}{(z^2+1)^2}A_3.$$
Therefore there are vector valued functions $A_4\equiv A_4(v)$ and $A_5\equiv A_5(v)$ depending only on $v$ such that:
$$A_3(v, z)=\frac{z}{\sqrt{z^2+1}}A_4(v)+\frac{1}{\sqrt{z^2+1}}A_5(v).\tag3.46 $$
The equation (3.41) becomes
$$x(u, v, z)=f(v, z)\cdot\Big(\cos(u)C_1+\sin(u)C_2\Big)+\frac{z}{\sqrt{z^2+1}}A_4(v)+\frac{1}{\sqrt{z^2+1}}A_5(v).$$
Differentiating the equation above and using (3.29), we have:
$$\align
0&=x_{vz}+\frac{z}{z^2+1}x_v\cr
&=\frac{1}{\sqrt{z^2+1}}A_4'(v).
\endalign$$
Therefore the vector valued function $A_4(v)$ is constant: $A_4(v)\equiv C_3\in\R^5$.

\medskip

\noindent
From the equation (3.42), one has
$$1=f^2+\frac{1}{z^2+1}\Big(z^2<C_3,C_3>+2z<C_3,A_5(v)>+<A_5(v),A_5(v)>\Big)$$
and thus
$$<C_3,C_3>=1, \quad <C_3,A_5(v)>=0,$$
and
$$<A_5(v),A_5(v)>=1-(1+z^2)f^2=1-\frac{e^{-2v}}{c_2}.\tag3.47 $$

\medskip

\noindent
Eliminating $y$ from the equations (3.32) and (3.34), we get
$$\align
0&=\alpha_1^2x_{vv}+f^{-2}x_{uu}+(z^2+1-\lambda^2)x_v-2z(z^2+1)x_z+2x\cr
&=\sqrt{z^2+1}\Big((c_2e^{2v}-1-c_1^2e^{4v})A_5''(v)+(1-c_1^2e^{4v})A_5'(v)+2A_5(v)\Big).
\endalign$$
Therefore the vector valued function $A_5(v)$ satisfies the following linear ordinary differential equation of second order:
$$(c_2e^{2v}-1-c_1^2e^{4v})A_5''(v)+(1-c_1^2e^{4v})A_5'(v)+2A_5(v)=0.\tag3.48$$
By the Lemma 3.4, one can conclude that the general solution of the equation (3.48) is
$$A_5(v)=g(v)C_4+h(v)C_5,$$
where $C_4, C_5\in\R^5$ are constant vectors, and $g$ and $h$ are the functions given in Remark 3.5.

\medskip

\noindent
But since from (3.47) the vector valued function $A_5$ satisfies $$<A_5(v), A_5(v)>=1-\frac{e^{-2v}}{c_2},$$ we have that $C_4$ and $C_5$ are orthonormal. They constitute together with $C_1, C_2, C_3$ an orthonormal basis of $\R^5$.
This proves Proposition 3.6. \hfill$\square$

\medskip

Our classifiaction theorem summarizes the results from Propositions 3.1-3.6.

\bigskip

{\bf Acknowledgements} The author would like to express his special thanks to
Professor An-Min Li and Professor Guosong
Zhao for their hospitality during his scientific visit in August 2001 at the
Sichuan University in Chengdu (China).

\bigskip

\bigskip

\noindent
\quad \quad \quad \quad \quad \quad \quad \quad \quad \quad \quad \quad \quad \quad \quad \quad \quad \quad \quad \quad \quad \quad \quad \quad Technische Universit\"at Berlin

\noindent
\quad \quad \quad \quad \quad \quad \quad \quad \quad \quad \quad \quad \quad \quad \quad \quad \quad \quad \quad \quad \quad \quad \quad \quad Institute f\"ur Mathematik 

\noindent
\quad \quad \quad \quad \quad \quad \quad \quad \quad \quad \quad \quad \quad \quad \quad \quad \quad
\quad \quad \quad \quad \quad \quad \quad Sekretariat MA 8-3 

\noindent
\quad \quad \quad \quad \quad \quad \quad \quad \quad \quad \quad \quad \quad \quad \quad \quad \quad \quad \quad \quad \quad \quad \quad \quad Stra\ss e des 17. Juni 136

\noindent
\quad \quad \quad \quad \quad \quad \quad \quad \quad \quad \quad \quad \quad \quad \quad \quad \quad \quad \quad \quad \quad \quad \quad \quad D-10623 Berlin, Germany

\medskip

\noindent
\quad \quad \quad \quad \quad \quad \quad \quad \quad \quad \quad \quad \quad \quad \quad \quad \quad \quad \quad \quad \quad \quad \quad \quad Email: lusala\@math.TU-Berlin.DE

\enddocument

\end